\documentclass[11pt]{article}

\title{Bernoulli-Taylor formula of $\psi$-umbral difference calculus}
\author{A.K.Kwa\'sniewski\\  
\\ Institute of Computer Science, Bia{\l}ystok University\\
PL-15-887 Bia{\l}ystok, ul.Sosnowa 64, POLAND
\\e-mail: kwandr@uwb.edu.pl}

\usepackage{amsmath,amsthm}

\chardef\bslash=`\\ 
\begin{document}
\maketitle

\begin{abstract}
We shall present here the $*_{\psi}$-Bernoulli-Taylor* formula of
a new sort  with the rest term of  the Cauchy type  recently
derived by the author in the case of $\psi$-difference calculus.
The central importance of such a type formulas is beyond any
doubt.
\end{abstract}
* see:  historical  remark at the beginning of this note
\\

\section{One Historical Remark} Here are the famous examples of expansion
$$\partial_0=\sum_{n=1}^{\infty}\frac{x^{n-1}}{n!}\frac{d^n}{dx^n}$$
or  $$\epsilon_0=\sum_{n=0}^{\infty}(-1)^n
\frac{x^{n}}{n!}\frac{d^n}{dx^n}$$ where $\partial_0$ is the
divided difference operator while $\epsilon_0$ is at the zero
point evaluation functional. If one compares these with
\textit{"series universalissima"} of J.Bernoulli from \textit{Acta
Erudicorum} (1694) (see commentaries in \cite{12}) and with
$$exp\{yD\}=\sum_{k=0}^{\infty}\frac{y^kD^k}{k!},\ \
D=\frac{d}{dx},$$ then confrontation with B.Taylor's {\em
"Methodus incrementorum directa et inversa"} (1715), London;
entitles one to call the expansion formulas considered in this
note {\em "Bernoulli - Taylor formulas"} or (for $n \rightarrow
\infty$) {\em "Bernoulli - Taylor series" \cite{1}}.

\section{Introduction}
While deriving  the Bernoulli-Taylor $\psi$-formula one is tempted
to adapt  the ingenious  Viskov`s method  [2] of arriving  to
formulas of such type  for various pairs of  operations. In our
case these would be $\psi$-differentiation and $\psi$-integration
(see: Appendix). However  straightforward application of  Viscov
methods in  $\psi$-extensions of umbral calculus leads to
sequences which are not normal (Ward) hence a new invention is
needed. This expected and verified here invention is the new
specific $*_{\psi}$ product of analytic functions or formal
series. This note is based on [3]  where the  derivation of this
new form of Bernoulli-Taylor $*_{\psi}$ - formula was delivered
due to the use of a specific $*_{\psi}$ product of formal series.

\section{Classical  Bernoulli-Taylor formulas with the rest term of  the Cauchy type by Viskov method}
Let us consider the obvious identity
\begin{equation} \label{eq1}
  \sum_{k=0}^{n} (\alpha_{k}-\alpha_{k+1})=\alpha_{0}-\alpha_{n+1}
\end{equation}
in which (\ref{eq1}) we now put
$\alpha_{k}=a^k b^k ; a,b\in \mathcal{A}.$
$\mathcal{A}$ is an associative algebra with unity over the field F=R,C.
Then we get
\begin{equation} \label{eq2}
\sum_{k=0}^n a^k (1-ab) b^k = 1-a^{n+1}b^{
n+1} ; a,b \in \mathcal{A}\
\end{equation}
Numerous choices of  $a,b \in \mathcal{A}$ result in many important
specifications of (\ref{eq2})

{\bf Example 1.} Let $\mathcal{F}$ denotes the linear space of sufficiently
smooth functions ${\it f: F \longrightarrow F}$. Let
\begin{eqnarray} \label{eq3}
a: \mathcal{F}\longrightarrow \mathcal{F};\quad &&(af)(x)= \int_{a}^{b}
f(t)dt,\nonumber \\
b: \mathcal{F}\longrightarrow \mathcal{F};\quad  &&(bf)(x) =
(\frac{d}{dx}f)(x);\\ l:\mathcal{F}\longrightarrow \mathcal{F};
\quad&& (lf)(x)=f(x)\nonumber.\end{eqnarray}
Then [b,a]=1-ab=$\varepsilon_\alpha$ where $\varepsilon_\alpha$ is evaluation
functional on $\mathcal{F}$ i.e.
\begin{equation}
\varepsilon_\alpha(f)=f(\alpha)
\end{equation}
\\
Using now the text-book integral Cauchy formula ($k>0$)
\begin{equation}
(a^k f)(x)= \int_a^x\frac{(x-t)^{k-1}}{(k-1)!} f(t)dt,\
\end{equation}
\\
and under the choice (\ref{eq3}) one gets from (\ref{eq2}) the well-known
Bernoulli-Taylor formula
\begin{equation} \label{eq6}
f(x)= \sum_{k=0}^n\frac{(x-\alpha)^k}{k!} f^{(k)} (\alpha)+R_{n+1}(x)
\end{equation}
\\
with the rest term $R_{n+1}(x)$ in the Cauchy form
\begin{equation} \label{eq7}
R_{n+1}(x)=\int_a^x\frac{(x-t)^{n}}{n!} f^{(n+1)}(t)dt
\end{equation}
\\
{\bf Example 2.} \cite{1} Let $\mathcal{F}$ denotes the linear space of
functions $f:Z_{+}\longrightarrow F; Z_{+}=N\cup \{0\}$. Let

\begin{eqnarray} \label{eq8}
a: Z_{+}\longrightarrow \mathcal{F};\quad &&(af)(x)= \sum_{k=0}^{x-1} f(k),
\nonumber \\
b: Z_{+}\longrightarrow \mathcal{F};\quad  &&(bf)(x) = f(x+1)-f(x),\\
l:Z_{+}\longrightarrow \mathcal{F};\quad&& (lf)(x)=f(x)\nonumber.
\end{eqnarray}
It is easy to see that  [b,a]=1-ab=$\varepsilon_0$ where
$\varepsilon_0$ is evaluation functional  i.e.
$\varepsilon_0(f)=f(0).$ $b=\Delta$ is the standard difference
operator with its left inverse definite summation operator a. The
corresponding $\Delta$ - calculus Cauchy formula is also known
(see formula (31 p.310 in \cite{5});
\begin{equation}
(a^k f)(x)= \sum_{r=0}^{x-1} \frac{(x-r-1)^{\underline {k-1}}}{(k-1)!}f(r);
k>0
\end{equation}
where $x^{\underline n}=x(x-1)(x-2)...(x-n+1).$
\\
Under the choice (\ref{eq8}) one gets from (\ref{eq2}) the $\Delta$ -
calculus Bernoulli - Taylor fomula \cite{1}
\begin{equation}
f(x)=\sum_{k=0}^n \frac {x^{\underline k}}{k!}(\Delta^k f)(0)+R_{n+1}(x)
\end{equation}
with the rest term $R_{n+1}(x)$ in the Cauchy $\Delta$ form
\begin{equation}
R_{n+1}(x)= \sum_{r=0}^{x-1} \frac{(x-r-1)^{\underline {n}}}{n!}
(\Delta^{n+1}f)(r);
\end{equation}

\section{$"\ast_\psi$ realization" of Bernoulli
identity.}  Now \textit{a specifically new} form of the
Bernoulli-Taylor formula with the rest term of the Cauchy type as
well as Bernoulli-Taylor series is to be supplied in the case of
$\psi$-difference \textit{umbral} calculus (see [5-8] and [9,10]
and references therein). For that to do we use natural
$\psi$-umbral representation [13,14] of Graves-Heisenberg-Weyl
(GHW) algebra  [11,12]  generators $\hat{p}$ and $\hat{q}$ and
then we use Bernoulli identity (\ref{eq12})
\begin{equation}  \label{eq12}
\hat p \sum_{k=0}^n \frac {(-\hat q)^k \hat p^k}{k!}=\frac{(-\hat q)^n \hat
p^{n+1}}{n!}
\end{equation}
derived by Viskov from (\ref{eq1}) under the substitution (see
(28) in [2])
$$\alpha_0=0,\ \alpha_k=(-1)^k{(\hat q)^{k-1} \hat p^k}{(k-1)!},\
k=1,2,...$$
due to $\hat p \hat q^n= \hat q^n \hat p +n \hat
q^{n-1}$  (n=1,2,...) resulting by induction from
\begin{equation} \label{eq13}
[\hat p, \hat q]=1
\end{equation}
\\{\bf Example 1.} The choice $\hat p=D \equiv \frac{d}{dx}$ and $\hat q
=\hat x-y, y\in F; \hat x f(x)=xf(x)$ after substitution into
Bernoulli identity (\ref{eq12}) and integration $\int_\alpha^x dt$
gives the Bernoulli - Taylor formula (\ref{eq6}).
\\{\bf Example 2.} The choice \cite{2} $\hat p=\Delta$  and $\hat q=
\hat x\circ E^{-1}$ where $E^{\alpha}f(x)= f(x+\alpha)$ after
substitution into Bernoulli identity (\ref{eq12}) and $"\Delta$ -
integration" $\sum_{r=0}^{\alpha-1}$ gives the Bernoulli - Mac
laurin formula of the following form ($\alpha, x \in {\bf Z},
\bigtriangledown =1-E^{-1})$ with the rest term $R_{n+1}(x)$
\begin{equation}
f(0)=\sum_{k=0}^{n} \frac{\alpha^{\underline{k}}}{k!}(-1)^{k+1}
(\bigtriangledown ^k f)(\alpha)+R_{n+1}(\alpha);
\end{equation}
\begin{equation}
R_{n+1}(\alpha)=(-1)^n\sum_{r=0}^{\alpha-1} \frac{r^{\underline n}}{n!}
(\bigtriangledown^{n+1}f)(r+1).
\end{equation}
\\{\bf Example 3.} Here  $f^{(k)}\equiv \partial_{\psi}^k f$ and $f(x)*_{\psi}g(x)\equiv f(\hat x_{\psi})g(x)$
 - see Appendix. The choice $\hat p =\partial_{\psi}$ and $ \hat q
= \hat z_{\psi}$  $(z=x-y)$ where $\hat x_{\psi}
x^n=\frac{n+1}{(n+1)_{\psi}}x^{n+1}$ after substitution into
Bernoulli identity (\ref{eq12}) and "$\partial_{\psi}$ -
integration" $\int_\alpha^x d_{\psi}t$ (see: Appendix) gives
another Bernoulli - Taylor $\psi$-formula of the form:
\begin{equation}
f(x)=\sum_{k=0}^n
\frac{1}{k!}(x-\alpha)^{k_{\ast_{\psi}}}\ast_{\psi}
f^{(k)}(\alpha)+ R_{n+1}(x)
\end{equation}
\\with the rest term $R_{n+1}(x)$ in the Cauchy-form
\begin{equation}
R_{n+1}(x)=\frac{1}{n!}\int_\alpha^x d_qt (x-t)^{n_{\ast_q}}\ast_q f^{(n+1)}
(t)dt
\end{equation}
\\In the above notation $x^{0*_{\psi}}=1,\ x^{n*_{\psi}}\equiv x*_{\psi} (x^{(n-1)*_{\psi}})=
x*_{\psi} ...*_{\psi} x=\frac{n!}{n_{\psi}!}x^n;\ n\geq 0.$
\\Naturally $\partial_{\psi} x^{n*_{\psi}}=nx^{(n-1)*_{\psi}}$ and in general $f,g$ - may be
formal series for which
\begin{equation}
\partial_{\psi}(f*_{\psi} g)=(Df)*_{\psi} g+f*_{\psi}(\partial_{\psi} g)
\end{equation}
\\i.e. Leibniz $*_{\psi}$ rule holds \cite{13,14,15}.
\\ Summary: These another forms of both the Bernoulli -Taylor
formula with the rest term of the Cauchy type \cite{3} as well as
Bernoulli - Taylor series are quite easily handy due to the
technique developed in \cite{13,14} where one may find more on
$*_\psi$ product devised perfectly suitable for the Ward's {\em
"calculus of sequences"} \cite{6} or more exactly $*_\psi$ is
devised perfectly suitable for the so-called $\psi$ - extension on
Finite Operator Calculus of Rota (see \cite{9,10,14,15} and
references therein)
\section{Appendix}  \textbf{$*_{\psi}$ product}\\
Let $n-{\psi}\equiv\psi_n$; $\psi_n\neq 0$: $n>0$.  Let
$\partial_\psi$ be a linear operator acting on formal series  and
defined accordingly by $\partial_\psi x^n=n_\psi x^{n-1}$. \\We
introduce now a intuition appealing
$\partial_\psi$-difference-ization rules for a specific new
$*_\psi$ product of functions or formal series. This $*_\psi$
product is what we  call:  the $\psi$-multiplication of functions
or formal series as specified below.
\\\textbf{Notation A.1.}\\
$x \ast _{\psi} x^{n} = \hat {x}_{\psi}(x^{n}) = \frac{{\left( {n
+ 1} \right)}}{{\left( {n + 1} \right)_{\psi} } }x^{n + 1};\quad n
\geq 0$ \; hence $x \ast _{\psi} 1 = (1_{\psi})^{-1}  x \not
\equiv x $ therefore $x  \ast_{\psi} \alpha 1 = \alpha 1 \ast
_{\psi}  x = x  \ast _{\psi} \alpha = \alpha \ast _{\psi}  x =
\alpha  (1_{\psi})^{-1} x$ and $\forall x, \alpha\in F$; $f(x)
\ast _{\psi} x^{n} = f(\hat {x}_{\psi})x^{n}$.
\\For $k \ne n $ \; x$^{n}
\ast _{\psi} $ x$^{k} \ne$ x$^{k} \ast _{\psi} $ x$^{n}$ as well
as x$^{n} \ast _{\psi} $ x$^{k} \ne$ x$^{n+ k}$ - in general.
\\In order to facilitate the  formulation  of observations  accounted
for on the basis of $\psi$-calculus representation of  GHW algebra
we shall use what follows. \\
 \textbf{Definition A.1.} With Notation A.1. adopted define the $*_\psi$
powers of $x$ according to
 $x^{n \ast _{\psi} } \equiv $ x $ \ast _{\psi} x^{\left( {n - 1} \right)
\ast _{\psi} } = \hat {x}_{\psi} (x^{\left( {n - 1} \right)\ast
_{\psi }} ) = $ x $ \ast _{\psi} $ x $ \ast _{\psi} $ ... $ \ast
_{\psi} $ x $=\frac{n!}{n_{\psi}  !}x^{n};\quad n \geq 0$. Note
that $x^{n\ast _{\psi} }  \ast _{\psi} x^{k\ast _{\psi} } =
\frac{{n!}}{{n_{\psi}  !}} x^{\left( {n + k} \right)\ast _{\psi} }
\ne x^{k\ast _{\psi} }  \ast _{\psi} x^{n\ast _{\psi} }  =
\frac{{k!}}{{k_{\psi}  !}}x^{\left( {n + k} \right)\ast _{\psi} }
$ for
$k  \ne n$ and $x^{0\ast _{\psi} }=1$.\\

This noncommutative $\psi$-product $ \ast _{\psi}$ is devised so
as to ensure the following observations.
\\\textbf{Observation A.1}
\begin{enumerate}
\renewcommand{\labelenumi}{\em \alph{enumi})}
\item $\partial _{\psi}  x^{n\ast _{\psi} } = n x^{\left( {n - 1}
 \right)\ast _{\psi} } $;\; $n \ge 0$
\item  {\bf exp}$_{\psi}  $[$\alpha ${\it x}] $ \equiv ${\bf exp}
\{$\alpha \hat {x}_{\psi}  $\}{\bf 1}
\item {\bf exp} [$\alpha x$] $ \ast _{\psi} $
({\bf exp}$_{\psi} $\{$\beta \hat {x}_{\psi}  $\}{\bf 1}) = ({\bf
exp}$_{\psi}  $\{[$\alpha +\beta $]$\hat {x}_{\psi}  $\}){\bf 1}
\item
  $\partial _{\psi} (x^{k} \ast _{\psi} \quad x^{n\ast _{\psi} } ) =
(D x^{k}) \ast _{\psi} x^{n \ast _{\psi} } + x^{k} \ast _{\psi}
(\partial _{\psi}  x^{n\ast _{\psi} })$

\item \label{e}
$\partial _{\psi} ( f \ast _{\psi} g) = ( Df) \ast _{\psi} g  + f
\ast _{\psi} (\partial _{\psi} g)$ ; $f,g$  - formal series
\item \label{f}
$f( \hat {x}_{\psi}) g (\hat {x}_{\psi} )$ {\bf 1} $= f(x) \ast
_{\psi} \tilde {g} (x)$ ; $\tilde {g} (x) = g(\hat
{x}_{\psi})${\bf 1}.
\end{enumerate}

\textbf{$\psi$-Integration} Let: $\partial_o x^n=x^{n-1}$. The
linear operator $\partial_o$ is identical with divided difference
operator. Let $\hat {Q}f(x)f(qx)$. Recall also that to the
"$\partial _{q} $ difference-ization" there corresponds the
$q$-integration which is a right inverse operation to
"$q$-difference-ization". Namely
\begin{equation}\label{eqa2}
 F\left( {z} \right): \equiv \left( {\int_{q} \varphi
}  \right)\left( {z} \right): = \left( {1 - q}
\right)z\sum\limits_{k = 0}^{\infty}  {\varphi \left( {q^{k}z}
\right)q^{k}}
\end{equation}
i.e.
\begin{multline}\label{eqa3}
 F\left( {z} \right) \equiv \left( {\int_{q} \varphi }
\right)\left( {z} \right) = \left( {1 - q} \right)z\left(
{\sum\limits_{k =
0}^{\infty}  {q^{k}\hat {Q}^{k}\varphi} }  \right)\left( {z} \right) =\\
=\left( {\left( {1 - q} \right)z\frac{{1}}{{1 - q\hat
{Q}}}\varphi} \right)\left( {z} \right).
\end{multline}
Of course
\begin{equation}
\partial _{q} \circ \int_{q} = id
\end{equation}
as
\begin{equation}
\label{eqa5} \frac{{1 - q\hat{Q}}}{{\left( {1 - q}
\right)}}\partial _{0} \left( {\left( {1 - q}\right)\hat
{z}\frac{{1}}{{1 - q\hat {Q}}}} \right)=id.
\end{equation}
Naturally (\ref{eqa5}) might serve to define a right inverse
operation to "$q$-difference-ization"
 $\left( {\partial _{q} \varphi}  \right)\left( {x} \right) = \frac{{1 - q\hat
{Q}}}{{\left( {1 - q} \right)}}\partial _{0} \varphi \left( {x}
\right)$ and consequently  the  "$q$-integration" as represented
by (\ref{eqa2}) and (\ref{eqa3}).  As it is well known  the
definite $q$-integral is an numerical approximation of  the
definite integral obtained in the $q \to 1$ limit.
\\Finally we introduce the analogous representation for $\partial_\psi$ difference-ization
\begin{equation}\label{eqa6}
  \partial_\psi=\hat n_\psi \partial_o;\ \hat n_\psi
  x^{n-1}=n_\psi x^{n-1};\ n\ge 1
\end{equation}
Then

\begin{equation}\label{eqa7}
  \int_\psi x^n=\left(\hat x \frac{1}{\hat
  n_\psi}\right) x^n=\frac{1}{(n+1)_\psi}x^{n+1};\ n\ge 0
\end{equation}
and of course  $\left(\int_\psi \equiv \int d_\psi \right)$
\begin{equation}\label{eqa8}
  \partial_\psi \circ \int_\psi=id
\end{equation}
Naturally $$ \partial_\psi \circ \int_a^x f(t)d_\psi t=f(x)$$

The formula of   "per partes"  $\psi$-integration  is easily
obtainable from (Observation A.1 e) and it  reads:

\begin{equation}\label{eqa9}
  \int_a^b(f*_\psi \partial_\psi g)(t)d_\psi t=[(f*_\psi
  g)(t)]_a^b - \int_a^b(Df*_\psi g)(t)d_\psi t
\end{equation}
\\
\textbf{Two Closing Remarks:}
\\
I.  All these above may be quite easily extended \cite {15} to the
case of any $Q\in End(P)$ linear operator that reduces by one the
degree of each polynomial \cite {16}. Namely one introduces \cite
{15}:

 \textbf{Definition A.2.}

    $$\hat {x}_{Q} \in End(P), \hat x_{Q}: F[x] \to F[x] $$
     such that  $(x^{n}) = \frac{{\left( {n
+ 1} \right)}}{{\left( {n + 1} \right)_{\psi} } }q_{n + 1}; n \geq
0;$ where $Qq_n=nq_{n-1}$.

Then   $\star_Q$ product of formal series  and  $Q$-integration
are defined analogously. (This has been accomplished by my student
E. Krot ).

II.   In  1937   Jean  Delsarte  [17]  had derived the general
Bernoulli-Taylor formula for a  class of linear operators $\delta$
including  linear operators that reduce by one the degree of each
polynomial. The rest term of the Cauchy-like type in his Taylor
formula (I) is given in terms of the unique solution of a first
order partial differential equation in two real variables. This
first order partial differential equation  is determined by the
choice of  the linear operator $\delta$  and the function  f under
expansion. In  our   Bernoulli-Taylor -formula (16)-(17) or in its
straightforward $\star_Q$     product of formal series and
$Q$-integration generalization -  there is no need to solve any
partial differential equation.

\end{document}